\documentclass[11pt]{article}

\usepackage{amscd, mathtools, amsmath, amssymb, fancyhdr, color}

\newcommand{\version}{version 2.0,\ \   5.04.2017}

\numberwithin{equation}{section}

\newcommand{\RR}{\mathbb{R}}
\newcommand{\ZZ}{\mathbb{Z}}

\newcommand{\la}{\lambda}
\newcommand{\al}{\alpha}
\newcommand{\be}{\beta}

\newcommand{\f}{\varphi}

\def\eqref#1{(\ref{#1})}

\newcommand{\arrow}{{\:\longrightarrow\:}}
\newcommand{\Z}{{\Bbb Z}}
\newcommand{\C}{{\Bbb C}}
\newcommand{\R}{{\Bbb R}}

\newcommand{\6}{\partial}
\def\1{\sqrt{-1}\:}
\newcommand{\restrict}[1]{{\left|_{{\phantom{|}\!\!}_{#1}}\right.}}
\newcommand{\cntrct}                
{\hspace{2pt}\raisebox{1pt}{\text{$\lrcorner$}}\hspace{2pt}}

\renewcommand{\tilde}{\widetilde}
\renewcommand{\bar}{\overline}
\renewcommand{\phi}{\varphi}
\renewcommand{\epsilon}{\varepsilon}
\renewcommand{\geq}{\geqslant}

\newcommand{\im}{\operatorname{im}}
\newcommand{\Lie}{\operatorname{Lie}}
\newcommand{\End}{\operatorname{End}}

\newcommand{\id}{\operatorname{\text{\sf id}}}

\newcommand{\Hom}{\operatorname{Hom}}

\newcommand{\const}{\operatorname{const}}

\newcommand{\coker}{\operatorname{coker}}

\newcommand{\Spec}{\operatorname{Spec}}

\def\ra{\arrow}

\def\m{\mathfrak m}
\def\N{\mathbb N}

\newcounter{Mycounter}[section]
\newcounter{lemma}[section]
\newcounter{claim}[section]
\newcounter{sublemma}[section]
\newcounter{corollary}[section]
\renewcommand{\thecorollary}{{Corollary
\thesection.\arabic{corollary}}}
\newcommand{\corollary}{%
     \setcounter{corollary}{\value{Mycounter}}
     \refstepcounter{corollary}
     \stepcounter{Mycounter}
     {\noindent \bf \thecorollary.\ }}

\newcounter{theorem}[section]
\renewcommand{\thetheorem}{{Theorem \thesection.\arabic{theorem}}}
\newcommand{\theorem}{%
     \setcounter{theorem}{\value{Mycounter}}
     \refstepcounter{theorem}
     \stepcounter{Mycounter}
     {\noindent \bf \thetheorem.\ }}

\newcounter{conjecture}[section]
\renewcommand{\theconjecture}{{Conjecture
\thesection.\arabic{conjecture}}}
\newcommand{\conjecture}{%
     \setcounter{conjecture}{\value{Mycounter}}
     \refstepcounter{conjecture}
     \stepcounter{Mycounter}
     {\noindent \bf \theconjecture.\ }}

\newcounter{proposition}[section]
\renewcommand{\theproposition}
       {{Proposition \thesection.\arabic{proposition}}}
\newcommand{\proposition}{%
     \setcounter{proposition}{\value{Mycounter}}
     \refstepcounter{proposition}
     \stepcounter{Mycounter}
     {\noindent \bf \theproposition.\ }}

\newcounter{definition}[section]
\renewcommand{\thedefinition}
       {{Definition~\thesection.\arabic{definition}}}
\newcommand{\definition}{%
     \setcounter{definition}{\value{Mycounter}}
     \refstepcounter{definition}
     \stepcounter{Mycounter}
     {\noindent \bf \thedefinition.\ }}

\newcounter{example}[section]
\newcounter{remark}[section]
\renewcommand{\theremark}{{Remark \thesection.\arabic{remark}}}
\newcommand{\remark}{%
     \setcounter{remark}{\value{Mycounter}}
     \refstepcounter{remark}
     \stepcounter{Mycounter}
     {\noindent \bf \theremark.\ }}

\newcounter{problem}[section]
\newcounter{question}[section]
\renewcommand{\leftmark}%
{{\scriptsize  Weighted Bott-Chern cohomology for LCK manifolds}}

\makeatletter

\@addtoreset{equation}{section}
\@addtoreset{footnote}{section}
\makeatother

\def\blacksquare{\hbox{\vrule width 5pt height 5pt depth 0pt}}
\def\endproof{\blacksquare}

\begin{document}
\begin{center}
{\LARGE\bf Weighted Bott-Chern and Dolbeault cohomology
 for LCK-manifolds with potential \\[3mm]
}

Liviu Ornea\footnote{\label{LO}Partially supported by CNCS UEFISCDI, project number PN-II
-ID-PCE-2011-3-
0118.}, Misha Verbitsky\footnote{Partially supported by RSCF grant 14-21-00053 within AG Laboratory NRU-HSE, and Simons-IUM fellowship.} and Victor Vuletescu\footnote{Partially supported by CNCS UEFISCDI, project number PN-II
-ID-PCE-2011-3-
0118.}

\end{center}

{\small
\hspace{0.15\linewidth}
\begin{minipage}[t]{0.7\linewidth}
{\bf Abstract} \\
A locally conformally K\"ahler (LCK) manifold is a 
complex manifold, with a K\"ahler structure on its
covering $\tilde M$, with the deck transform group acting
on $\tilde M$ by holomorphic homotheties.
One could think of an LCK manifold as of a 
complex manifold with a K\"ahler form taking
values in a local system $L$, called {\bf the conformal 
weight bundle}. The $L$-valued cohomology of $M$
is called {\bf Morse-Novikov cohomology}; it was
conjectured that (just as it happens for K\"ahler
manifolds) the Morse-Novikov complex satisfies
the $dd^c$-lemma, which (if true) would have far-reaching
consequences for the geometry of LCK manifolds. 
In particular, this version of $dd^c$-lemma would 
imply existence of LCK potential on any LCK manifold with
vanishing Morse-Novikov class of its 
$L$-valued Hermitian symplectic form. The 
$dd^c$-conjecture was disproved for Vaisman manifolds 
by R. Goto. We prove that $dd^c$-lemma is true
with coefficients in a sufficiently general power of $L$
on any Vaisman manifold or LCK with potential.  
\end{minipage}
}

\tableofcontents

\section{Introduction}


\subsection{LCK manifolds and $d_\theta d_\theta^c$-lemma}

A locally conformally K\"ahler (LCK) manifold is a 
complex manifold which admits a K\"ahler metric 
on its universal covering $\tilde M$ such that the monodromy
acts on $\tilde M$ by K\"ahler homotheties. For more
details and the reference on this subject, please see
Section \ref{_LCK_details_Section_}. 

The LCK property is equivalent to existence
of a Hermitian form $\omega$ on $M$ satisfying
$d\omega=\omega\wedge \theta$, where $\theta$
is a closed 1-form. This form is called
{\bf the Lee form} of an LCK-manifold.

One can consider
the K\"ahler form on $\tilde M$ as a K\"ahler form on
$M$ taking values in a 1-dimensional local system,
or, equivalently, in a flat line bundle $L$. This bundle
is called {\bf the weight bundle} of $M$.

The cohomology of this local system is known as 
{\bf the Morse-Novikov cohomology} of an LCK manifold.
In locally conformally K\"ahler geometry, the 
Morse-Novikov cohomology shares many properties 
of the Hodge decomposition with
the usual cohomology of the complex manifolds.
The locally conformally K\"ahler form represents
a cohomology class (called the Morse-Novikov class)
of an LCK manifold, encoding the topological
properties of an LCK structure. However, the 
$dd^c$-lemma, which plays a crucial role 
for the K\"ahler geometry, is invalid in 
the Morse-Novikov setting. The main question 
of the locally conformally K\"ahler geometry
is to find a replacement of the $dd^c$-lemma
which would allow one to study the interaction
between the complex geometry and the topology of a manifold.

The statement of the $dd^c$-lemma seems,
on the first sight, to be technical.
It says that on any compact K\"ahler manifold
$(M,I)$, one has $\im d \cap \ker d^c = \im dd^c$,
where $d^c=IdI^{-1}$ is the twisted de Rham differential.
However, it is used as a crucial step in the proof
of the degeneration of the Dolbeault-Fr\"olicher
spectral sequence, and in the proof of 
homotopy formality of K\"ahler manifolds.

For an LCK manifold, one replaces the de Rham differential
by its Morse-Novikov counterpart $d_\theta:=d-\theta$,
where $\theta$ is the connection form of its
weight bundle; the twisted de Rham differential
is replaced by $d^c_\theta=Id_\theta I^{-1}$.
It was conjectured in \cite{_OV:MN_} that
the $d_\theta d^c_\theta$-lemma would hold on 
any LCK manifold, giving $\im d_\theta \cap \ker
d^c_\theta = \im d_\theta d^c_\theta$.
The implication of the $d_\theta d^c_\theta$-lemma
would include the topological classification of
LCK structures on some manifolds (such as nilmanifolds) 
and a construction of automorphic K\"ahler potentials
on LCK manifolds with vanishing Morse-Novikov class.
However, this conjecture was false, as shown by 
R. Goto (\cite{goto}).

\subsection{Weighted Bott-Chern cohomology}

When the $d_\theta d^c_\theta$-lemma is false, one needs
to study a more delicate cohomological invariant, 
called {\bf the weighted Bott-Chern cohomology of a manifold:}
\[
H^{p,q}_{BC}(M,L):= 
\frac{\ker d_\theta \cap \ker d^c_\theta}
{\im d_\theta d^c_\theta}\restrict{\Lambda^{p,q}(M)}.
\]
In \cite{goto}, Goto has shown that 
the Bott-Chern cohomology group is responsible
for the deformational properties of an LCK manifold,
and computed it for certain $(p,q)$ and certain
examples of LCK manifolds, called {\bf the Vaisman
manifolds} (see Subsection \ref{_Vaisman_Subsection_}).

\hfill

\definition\label{_weight_power_Definition_}
The local system $L$ associated to a LCK manifold $M$ is a real, oriented line bundle over $M$ 
with a flat connection. Trivializing this bundle,
we can write its connection as $\nabla_L= d -\theta$,
where $\theta$ is the Lee form of our LCK manifold.
For arbitrary $a\in \C$, the connection 
 $\nabla_{L_a}:= d-a\theta$
is also flat. For $a\in \Z$, the corresponding
line bundle is identified with the $a$-th tensor power of $L$,
denoted as $L^a$. One may think of the flat line
bundle $(L, \nabla_{L_a})$ as of a real (or complex)
power of $L$.  We denote this line bundle and its
local system by $L_a$, and call it
{\bf $a$-th power of a weight bundle.}

\hfill

In this paper we compute the weighted Bott-Chern cohomology 
for $L_a$, on LCK manifolds with proper potential,
and show that it vanishes for all $a$
outside of a discrete countable subset of $\R$ (\ref{_BC_vanishing_Corollary_}).
This implies $dd^c$-lemma for forms with coefficients
in $L_a$, for these values of $a$. This result is based on a computation of 
Dolbeault cohomology with coefficients in $L_a$, which
also vanishes for all $a$ but a discrete countable subset (\ref{mainresult}).

\subsection{LCK manifolds with potential}

\definition
A compact LCK manifold $(M,\omega, \theta)$ 
is called {\bf LCK with potential} if $\omega= d_\theta
d^c_\theta \psi$ for a positive function $\psi$ which is called
{\bf LCK potential}. 

\hfill

An equivalent definition will be given in Subsection \ref{lpot}.

LCK manifolds with potential are understood very well now.
The following results were proven in
\cite{_OV:_Potential_} and \cite{_OV:_Topology_imrn_} (see also \cite{ov_rank}).
Recall that {\bf a linear Hopf manifold} is a quotient
of $\C^n \backslash 0$ by a $\Z$-action generated
by a linear map with all eigenvalues $|\alpha_i| >1$.

\hfill

\theorem
Let $M$ be a compact complex manifold.
Then $M$ admits an LCK metric with potential
if and only if $M$ admits an embedding to
a linear Hopf manifold.
\endproof

\hfill

\theorem
Let $M$ be an LCK manifold with potential.
Then $M$ is a deformation of a Vaisman manifold
(\ref{_Vaisman_Definition_}). In particular,
$M$ is diffeomorphic to a principal $S^1\times S^1$-bundle
over a projective orbifold.
\endproof

\hfill

It would be nice to have a topological characterization
of LCK manifolds with potential. Since \cite{_OV:MN_},
we were extending much effort trying to prove the
following conjecture, which has  many geometric
consequences.

\hfill

\conjecture\label{_MN=0_then_potential_Conjecture_}
Let $(M,\omega,\theta)$ be a compact LCK manifold. 
Assume that $\omega$ is $d_\theta$-exact.
Then $\omega$ is $d_\theta d^c_\theta$-exact,
that is, $M$ is a LCK with potential.

\hfill

This conjecture is still open. It would trivially
follow if the $d_\theta d^c_\theta$-lemma were true,
but it is known now to be false. However, a weaker conjecture
still stands.

\hfill

\conjecture\label{_dd^c_generic_a_Conjecture_}
Let $(M,\omega,\theta)$ be a compact LCK manifold,
$L$ its weight bundle, and $L_a$ the weight bundle
to the power of $a\in \R$ (\ref{_weight_power_Definition_}).
Then, for all $a$ outside of a discrete countable
set, $d_{a\theta} d^c_{a\theta}$-lemma is true:
for any $d_{a\theta}$-exact (1,1)-form $\eta$,
one has $\eta= d_{a\theta} d^c_{a\theta}f$ (but this does not imply that the $d_{a\theta} d^c_{a\theta}$-lemma is true for other bidegrees).

\hfill

In this paper, we prove that
\ref{_dd^c_generic_a_Conjecture_}
is true for LCK manifolds with proper potential
(\ref{_BC_vanishing_Corollary_}). This is done by first proving a generic vanishing result for weighted Dolbeault cohomology (\ref{mainresult}).


\section{Locally conformally K\"ahler geometry}
\label{_LCK_details_Section_}


In this section we give the necessary definitions and
properties of locally conformally K\"ahler (LCK)
manifolds.

\subsection{LCK manifolds}

\definition
A complex manifold $(M,I)$ is LCK if it admits
a K\"ahler covering $(\tilde M, \tilde \omega)$,
such that the covering group acts by holomorphic
homotheties.

\smallskip

Equivalently, there exists on $M$ a \emph{closed} 1-form
$\theta$, called \emph{the Lee form}, such that $\omega$
satisfies the integrability condition:
$$d\omega=\theta\wedge\omega.$$

Clearly, the metric $g:=\omega(\cdot, I\cdot)$ on $M$ is
locally conformal to some K\"ahler metrics and its lift to
the K\"ahler cover in the definition is globally conformal
to the K\"ahler metric corresponding to $\tilde\omega$.

To an LCK manifold one associates the \emph{weight bundle}
$L_{\mathbb{R}}\longrightarrow M$.
It is a real line bundle associated to the
representation\footnote{ In conformal geometry, the weight
bundle usually corresponds to $\mid\det
A\mid^{\frac{1}{2n}}$. For LCK-geometry, $\mid\det A\mid^{\frac{1}{n}}$
is much more convenient.}
\[
\mathrm{GL}(2n,\RR)\ni A\mapsto \mid \det A\mid^{\frac{1}{n}}.
\]
The Lee
form induces a connection in $L_{\mathbb{R}}$ by the
formula $\nabla=d- \theta$. $\nabla$  is
associated to the Weyl covariant derivative (also denoted
$\nabla$) determined on $M$ by the
LCK metric and the Lee form.
As $d\theta=0$, then $\nabla^2= d\theta=0$, and hence
$L_{\mathbb{R}}$ is flat.

The complexification of the weight bundle will be denoted by $L$. The Weyl
connection extends naturally to $L$ and its
$(0,1)$-part endows $L$ with a holomorphic structure.

\subsection{Vaisman manifolds}
\label{_Vaisman_Subsection_}

\definition\label{_Vaisman_Definition_}
A Vaisman manifold is an LCK manifold with $\nabla^g$-parallel Lee form,
where $\nabla^g$ is the Levi-Civita connection.

\hfill

The following definition is implicit in the work of Boyer \& Galicki, see \cite{bg}:

\hfill

\definition
A Sasakian manifold is an odd-dimensional contact manifold
$S$ such that its symplectic cone $CS$ is equipped with
a K\"ahler structure, compatible with its symplectic
structure, and the standard symplectic homothety 
map $\rho_t:\; CS \arrow CS$ is holomorphic.

\hfill

Compact Vaisman manifolds can be described in terms of
Sasakian geometry as follows.

\hfill

\theorem
Let $(M,I,g)$ be a compact
Vaisman manifold. Then 
$M$ admits a conic K\"ahler covering
$(W\times \mathbb{R}_+, t^2g_W+dt^2)$ such that  the covering group is
an infinite cyclic group, generated by the
transformation $(w,t)\mapsto (\varphi(w), qt)$
for some Sasakian automorphism $\varphi$ and $q\in\mathbb{Z}$.

\hfill

The typical example of a compact Vaisman manifold is the
diagonal Hopf manifold
$H_A:=\C^n/ \langle A\rangle$ with
$A=\text{diag}(\al_i)$, with $|\alpha_i|<1$. An explicit construction of the Vaisman metric 
on $H_A$ is given in \cite{ov_shells}. Other Vaisman metrics appear on compact complex surfaces, \cite{_Belgun_}. 

Among the LCK manifolds which 
do not admit Vaisman metrics are some of the Inoue surfaces (cf.
\cite{_Tricerri_}, \cite{_Belgun_}) and their generalizations to higher
dimensions (\cite{_Oeljeklaus_Toma_}). The rank 0 Hopf surfaces
are also non-Vaisman (\cite{_Gauduchon_Ornea_}).

\subsection{LCK manifolds with potential}\label{lpot}

\definition \label{_LCK_w_pote_Definition_}
(\cite{_OV:_Potential_})
A compact complex manifold $(M,I)$ is {\bf LCK  with  potential} if it admits a
K\"ahler cover $(\tilde M, \tilde \omega)$ with global potential $\f:\; \tilde M
\rightarrow \RR_{+}$,
such that  and the monodromy  map $\tau$ acts on $\f$ by multiplication with a
constant: $\tau (\f)=\const \cdot \f$.

If $\f$ is {\bf proper} (inverse images of compact sets are compact), then $(M,I)$ is called {\bf LCK with proper potential}.

\hfill

\remark\label{_OV_defi_Remark_}

In 
\cite[Proposition 2.5]{_OV:_Potential_} (see also \cite{ov_rank}) it was proven that
$\f$ is proper if and only if the monodromy of the weight
bundle is discrete in $\R_{+}$, that is, isomorphic to 
$\Z$. 

\hfill

Vaisman manifolds are LCK with potential 
(the potential is equal to the squared norm of the Lee
field), which can be easily seen from the Sasakian
description given above (\cite{_Verbitsky_vanishing_}).
LCK metrics with potential are in one to one
correspondence with strongly pseudoconvex shells in 
affine cones, as shown in \cite{ov_shells}.

\hfill

We summarize the main properties of compact LCK manifolds with potential:

\hfill

\theorem \label{pot}
\begin{description}

\item[(i)] (\cite{_OV:_Potential_}) The class of compact LCK manifolds with
potential is stable to small deformations. 

\item[(ii)] 
(\cite[Theorem 2.1]{_OV:_Topology_imrn_})
Any LCK manifold with potential can be deformed to
a Vaisman manifold. Moreover, the
set of points which correspond to Vaisman manifolds is dense in the moduli of 
compact LCK manifolds with potential 

\item[(iii)] (\cite{_OV:_Potential_}) Any compact LCK manifold with potential
 can be holomorphically embedded into a
Hopf manifold.  Moreover, a compact Vaisman manifold 
can be holomorphically embedded in a diagonal
Hopf manifold.
\end{description}

\subsection{Morse--Novikov complex and cohomology of local systems}
\label{_MN_Subsection_}

Let $M$ be a smooth manifold, and $\theta$ a closed 1-form
on $M$. Denote by $d_\theta:\; \Lambda^i(M)\arrow \Lambda^{i+1}(M)$
the map $d-\theta$. Since $d\theta=0$, $d_\theta^2=0$.

Consider the {\bf  the Morse--Novikov complex},
(see {\em e.g.}   \cite{_Pajitnov_}, \cite{_Ranicki_},
\cite{_Millionschikov_})
\[
\Lambda^0(M)\stackrel {d_\theta} \arrow \Lambda^{1}(M)
\stackrel {d_\theta} \arrow\Lambda^{2}(M) \stackrel {d_\theta} \arrow \cdots
\]
Its cohomology is {\bf the Morse--Novikov cohomology} of $(M,\theta)$.

In Jacobi and locally conformal symplectic  geometry, this object is called
{\bf Lichnerowicz-Jacobi}, or {\bf Lichnerowicz cohomology}, motivated by
Lichnerowicz's work \cite{_Lichnerowicz_} on Jacobi manifolds (see {\em e.g.}  \cite{_Leon_Lopez_} and \cite{Banyaga}).

\hfill

Obviously, the flat line bundle $L$ can be viewed as a local system associated with the character $\chi:\pi_1(M) \arrow \R^{>0}$
given by the exponential $e^\theta \in H^1(M, \R^{>0})$,
considered as an element of $\R^{>0}$-valued cohomology. Then we have:

\hfill

\proposition\label{loc_syst=mn}
(see {\em e.g.}  \cite{Novikov})
The cohomology of the local system $L$
is naturally identified with the
cohomology of the Morse--Novikov
complex $(\Lambda^*(M), d_\theta)$.

\hfill

The following result was proven in \cite{_Leon_Lopez_} and, with a different method, in \cite{_OV:MN_}:

\hfill

\theorem\label{MN} The Morse--Novikov cohomology of a compact Vaisman manifold vanishes identically.

\hfill

On the other hand on one of the Inoue surfaces (which is
LCK but non-Vaisman) the Morse--Novikov class of $\omega$
is non--zero, see \cite[Theorem 1]{Banyaga}.

\section{Weighted Dolbeault cohomology for LCK manifolds with potential}


Let $M$ be an LCK manifold with proper potential, and $\tilde M$
its $\Z$-covering equipped with the automorphic K\"ahler metric.
In \cite{_OV:_Potential_} it was shown that the metric 
completion $\tilde M_c$ of $M$ is a Stein variety with 
at most one isolated singularity. Moreover, 
$\tilde M_c$ is obtained from $\tilde M$ by adding
one point, called ``the origin''. Denote this point by $c$.

\hfill

\remark 
If $M$ is Vaisman, then $\tilde M$ is a true 
(Riemannian) cone and the fibres are Sasakian. In the general case, nothing 
more precise can be said neither on the metric of 
$\tilde M$ nor on the contact metric structure of the fibres.

\hfill

 Since $\tilde M_c$ a singular variety, to control what happens in the neighbourhood of $c$ we need 
 some technique borrowed from algebraic geometry which we briefly explain below. Note that we could arrive at the 
 same results by using $L^2$-estimates, but the computations and technicalities would have been much more involved.

\subsection{Main result: the generic vanishing theorem}
The main result of this paper is:

\theorem\label{mainresult}
Let $M$ be an LCK manifold with proper potential, $\theta$ its Lee form, $\tilde M$ its K\"ahler  $\Z-$cover and denote by $t:\tilde M\ra \tilde M$ the monodromy action. 
Let  $\alpha\in \C$ be arbitrary and let $L_\alpha$ be the flat line bundle on $M$ corresponding to $\alpha\!\cdot \!\theta.$

Then  for any $q\in \N$ 
$$H^q(M, \Omega^p_M\otimes L_\alpha)=0,$$  for all $\alpha\in \C$ but a discrete countable subset.

\hfill

\remark
For  some Hopf manifolds,  stronger  vanishing results  were obtained by Ise \cite{_Ise_} and Mall \cite{_Mall_}.  In these cases, the set of exceptions is made explicit. 

\hfill

We  describe the main steps of the proof and give the details in the next section.

\hfill

\noindent {\bf Step 1: reduction to the local cohomology.} One has the following exact sequence,  (see \ref{_exa_long_equi_forms_Corollary_}, which follows from \ref{gode}):
\begin{multline}\label{exa1}
0\arrow H^0(M,\Omega^i_M \otimes L_\alpha) 
\arrow H^0(\tilde M,\Omega^i_{\tilde M})\stackrel{t-\alpha}
\arrow \\
\stackrel{t-\alpha}
\arrow H^0(\tilde M,\Omega^i_{\tilde M})\arrow\allowbreak
H^1(M,\Omega^i_M \otimes L_\alpha)
\arrow \cdots
\end{multline}
We are thus reduced to the study of the maps 
$$H^j\left(\tilde M,\Omega^i_{\tilde M}\right)\stackrel{t-\alpha}
\arrow 
H^j\left(\tilde M, \Omega^i_{\tilde M}\right).$$

Denote by  $\Omega^i_{\tilde M_c}$ be the exterior $i-$power of the sheaf of K\"ahler differentials on $\tilde M_c$ and {by} $S$ its stalk at $c.$ Using cohomology with supports, we have an exact sequence
\begin{multline*}
0\arrow H^0_{\m}(S) 
\arrow H^0\left(\tilde M_c, \Omega^i_{\tilde M_c}\right)\arrow
H^0\left(\tilde M, \Omega^i_{\tilde M}\right)\arrow\\
\arrow
H^1_{\m}(S)
\arrow H^1\left(\tilde M_c, \Omega^i_{\tilde M_c}\right)\arrow\cdots
\end{multline*}
Since $\tilde M_c$ is Stein, $H^j\left(\tilde M_c, \Omega^i_{\tilde M_c}\right)=0$ for all $j\geq 1$,   we obtain isomorphisms 
$$H^j\left(\tilde M, \Omega^i_{\tilde M}\right)\simeq H^{j+1}_\m(S),$$  
and an exact sequence
\begin{multline*}
0\arrow H^0_{\m}(S) 
\arrow H^0\left(\tilde M_c, \Omega^i_{\tilde M_c}\right)\stackrel{t-\alpha}\arrow
H^0\left(\tilde M, \Omega^i_{\tilde M}\right)\arrow
H^1_{\m}(S)
\arrow 0
\end{multline*}

These induce the commutative diagrams
\begin{equation}\label{diag1}
\begin{CD}
H^i\left(\tilde M, \Omega^j_{\tilde M}\right)@>{\simeq}>> H^{i+1}_\m(S)\\
@VVV @VVV\\
H^i\left(\tilde M, \Omega^j_{\tilde M}\right)@>{\simeq}>> H^{i+1}_\m(S)
\end{CD}
\end{equation}
and respectively
{\scriptsize \begin{equation}\label{diag2}
\begin{CD}
0@>>>H^0_{\m}(S)@>>>H^0\left(\tilde M_c, \Omega^i_{\tilde M_c}\right)@>>> H^0\left(\tilde M, \Omega^i_{\tilde M}\right)@>>>H^1_{\m}(S)@>>>0\\
@. @V{t-\alpha}VV @V{t-\alpha}VV @V{t-\alpha}VV @V{t-\alpha}VV @.\\
0@>>>H^0_{\m}(S)@>>>H^0\left(\tilde M_c, \Omega^i_{\tilde M_c}\right)@>>> H^0\left(\tilde M, \Omega^i_{\tilde M}\right)@>>>H^1_{\m}(S)@>>>0
\end{CD}
\end{equation}}

Eventually, notice that $H^{i+1}_{\m}(S)$ and $H^0\left(\tilde M_c, \Omega^i_{\tilde M_c}\right)$ are $R-$modules.

\hfill

\noindent {\bf Step 2: algebraic proof of generic vanishing.} At this step we use the following result, which will be proven in section \ref{subsec:AlgProof}:

\hfill

\theorem\label{gene_iso} For any local Noetherian $\C-$algebra $R$ endowed with a $\Z-$action given by an automorphism of local $\C-$algebras $ t_R$ and for any $R-$module $N$ endowed also with a $\Z$ action $t_N$ which is $ t_R$-equivariant, {\em i.e.}  
$$t_N(rm)= t_R(r) t_N(m), \quad \text{for all}\,\, r\in R ,\, m\in N,$$  
the map $t_M-\alpha$ is a $\C-$linear isomorphism for all $\alpha \in \C$ but a countable subset.

\hfill

\noindent{\bf Step 3.} Using the above commutative diagrams \eqref{diag1}, \eqref{diag2}, we conclude that for each $\alpha \in \C$ but a countable subset and any $i, j\geq 0$ the map
$$t-\alpha:
H^i\left(\tilde M, \Omega^j_{\tilde M}\right)
\arrow
H^i\left(\tilde M, \Omega^j_{\tilde M}\right)
$$
is an isomorphism. From the exact sequence (\ref{exa1}) we obtain $H^i(M,\Omega^j_M\otimes L_\al)=0$, for all $\al$ in $\C$ but a countable set. Moreover, by upper-continuity on $\al$, the set $\{\al\in\C\,;\, H^i(M,\Omega^j_M\otimes L_\al)=0\}$ is analytically Zariski open, and hence its complement is discrete since it is countable.
\hfill \endproof

\subsection{Proof of Step 1: reduction to the local cohomology}

\newcommand{\God}{\operatorname{God}}
\definition
Let $F$ be a sheaf of $\C$-vector spaces over a topological vector space.
Denote by $F_x$ the stalk of $F$ in $x\in M$, and 
let $\God(F)$ be the sheaf defined by 
$\God(F)(U):= \prod_{x\in U}F_x$.
The natural sheaf embedding $F \hookrightarrow \God(F)$
is apparent. The sheaves $\God_i(F)$ are defined
inductively: 
set $\God_0(F):=F,\, \God_1(F):=\God(F)$, and then   
$$\God_{i+1}(F):= \God(\God_i(F)/\God_{i-1}(F)).$$ 
This gives an exact sequence
\[
0\arrow F \arrow \God_1(F)\arrow \God_2(F)\arrow\cdots
\]
called {\bf the Godement resolution of $F$}.

\hfill

\theorem\label{godement_resolution}\label{gode}
Let $\tilde M\stackrel \pi \arrow M$ 
be a manifold equipped with a free action of $\Z$,
$M:= \tilde M/\Z$ its quotient, and let $F$ be a
$\Z$-equivariant sheaf on
$\tilde M$. For any character $\alpha:\; \Z \arrow \R$, 
denote by $F_\alpha\subset \pi_* F$ the sheaf of automorphic sections
of $\pi_*F$, associated with the character $\alpha$, considered
as a sheaf on $M$. 

Then one has the exact sequence
\begin{equation}\label{_long_exact_equiva_Equation_}
0\arrow H^0(M,F_\alpha) \arrow H^0(\tilde M, F)\stackrel{t-\alpha}
\arrow  H^0(\tilde M, F)\arrow H^1(M,F_\alpha)\arrow \cdots
\end{equation}
where $t$ is the associated action by 
the generator of $\Z$ acting on $\tilde M$,
and $\alpha$ is the multiplication by the number 
$\alpha(t)$.

\hfill

\noindent{\bf Proof:} 
Consider the Godement resolution
$0 \arrow F \arrow F^1 \arrow F^2 \arrow \cdots$. Here $F^i=\mathrm{God}(F^{i-1}/\mathrm{im}(d_{i-1}))=\mathrm{God}(\coker(d_{i-1}))$, $F^0=F$, and $d_i:F^{i-1}\arrow F^i$.  
Then 
\begin{equation}\label{_short_exact_gode_}
0 \arrow F^k_\alpha \arrow \pi_* F^k \stackrel{t-\alpha}
\arrow \pi_k F^*\arrow 0
\end{equation}
is an exact sequence of complexes of flabby sheaves over
$M$. 

Indeed, $F^k_\alpha=\ker (t-\al)$ and we only have to show
that $t-\al$ is surjective. It is enough to make the proof
at the level of sections of $F^k$. 
The argument is combinatorial. We look at $\tilde M$ as
$\bigcup_{i\in\Z}\tilde M_i$ where $M_0$ is a 
fundamental domain of the $\Z$ action and $\tilde
M_i=t^i(M_0)$. 

Then, given $f\in F^k(U)$, $U\subset \tilde M$, it is enough
to solve the equation $(t-\al)g=f$ for each
$f_i=f|_{U_i}$, $U_i=U\cap \tilde M_i$; this will 
give as solution the section $g_{i-1}\in F(U_{i-1})$,
$i\in\Z$. The equation is 
$$tg_it^{-1}-\al g_{i-1}=f_{i-1},$$
which can be solved recursively once we have chosen arbitrarily $g_0\in F(U_0)$.
 
The  long exact sequence associated to \eqref{_short_exact_gode_} is precisely
\eqref{_long_exact_equiva_Equation_}.
\hfill\endproof

\hfill

Let now $M$ be a locally conformally K\"ahler manifold with 
 K\"ahler covering $\tilde M$ and monodromy $\Gamma\cong \Z$.
Consider the weight bundle $L$ on $M$, and let $L_\alpha$ be its
power associated with the character $\alpha \in \Hom(\Gamma, \R)$. 
Since the automorphic forms on 
$\tilde M$ can be identified with forms on $M$ 
with values in $L$, from the above result we directly obtain:

\hfill

\corollary\label{_exa_long_equi_forms_Corollary_} 
For a compact LCK manifold with monodromy $\Z$ one has the exact sequence 
for the Dolbeault cohomology of $M$ with values in $L_\alpha$:
\begin{multline*}
0\arrow H^0(M,\Omega^i_M \otimes L_\alpha) 
\arrow H^0\left(\tilde M, \Omega^i_{\tilde M}\right)\stackrel{t-\alpha}
\arrow \\ \stackrel{t-\alpha}
\arrow H^0\left(\tilde M, \Omega^i_{\tilde M}\right)\arrow H^1(M,\Omega^i_M \otimes L_\alpha)
\arrow \cdots
\end{multline*}


\subsection{Proof of Step 2: algebraic proof of generic vanishing.} \label{subsec:AlgProof}

\begin{remark}\label{pix}
Let  $(V_n, t_n)_{n\geq 0}$ be a sequence of finite-dimensional vector spaces and endomorphisms $t_n:V_n\ra V_n.$
Let $V=\prod_{n\geq 0} V_n$ and $t=\prod_{n\geq 0} t_n.$ Then 
$$\Spec(t)=\bigcup_{n\geq 0} \Spec(t_n)$$
In particular, $\Spec(t)$ is at most countable.

Here, for  a $\C$-vector space $V$ and $u\in\End(V)$, $\Spec(u):=\{\la\in\C$ \,;\,  $u-\la\cdot\id\, \text{ is not an isomorphism}\}$.
\end{remark}

\hfill

This implies the following:

\hfill

\begin{lemma}\label{rema} If  $(M, t_m)$ is a finitely generated complete $R-$module which is equivariant,
then $\Spec(t_M)$ is at most countable.
\end{lemma}

\hfill

\noindent {\bf Proof:}  Since $M$ is complete we have
$$M=\prod_{n\geq 0}\m^nM/\m^{n+1}M.$$ Since $M$ is finitely generated, $\m^nM/\m^{n+1}M$ is finite dimensional $\C$-vector space  for all $n\geq 0$, so \ref{pix} applies.
\hfill\endproof

\hfill

Unfortunately, the cohomology modules $H^i_{\m}(M)$ are usually not finitely generated, so we need to elaborate further, by first reducing to the case of regular rings, and then using local duality and the explicit description of the injective hull of the residue field.

First, since local cohomology does not change under completion (cf \cite{Hun}, Prop. 2.15),  we may assume that both $R$ and $M$ are complete. 

Next, we reduce  to the case when $R$ is regular.

To do this, we choose a minimal system of generators for $\m_R$,  $m_1,\ldots, m_n$ and define a map
$$\pi:S=\C[[X_1,\dots, X_n]]\ra R,$$ 
by $X_i\mapsto m_i, i=1,n.$

The action $t_R$ on $R$  lifts to  an action $t_S$ on $S$ as follows. Choose lifts $s_i\in S$ of $t(m_i)$ for all $i=1,\dots,n$, and define $t_S(X_i)=s_i$. Note that $t_S$ is well-defined as a morphism of local $\C$-algebras by \cite[Theorem 7.16]{eis}.

So we can look at $M$ as an equivariant $S$-module.

Also, the local cohomology is preserved, 
since $\m_R=\m_SR $ and  using \cite[Proposition 2.14 (2)]{Hun},  we have
$H^i_{\m_S}(M)\simeq H^i_{\m_R}(M).$

Denote by $t^i$ the endomorphism of $H^i_\m(M)$ induced by $t_M$ and $t_R$.

By local duality (\cite[Theorem 4.4]{Hun}) we have:
$$H^i_{\m}(M)\simeq \mathrm{Ext}^{n-i}_R(M, R)^{\vee}=\mathrm{Hom}_R(\mathrm{Ext}^{n-i}_R(M, R), E(k))$$
where $E(k)$ is the injective hull of the residue field.

For regular rings, the injective hull $E(k)$ is described by Lyubeznik (\cite{_liu_}):
$$E(k)={\mathcal D}/\mathfrak m{\mathcal D}$$
where ${\mathcal D}$ is the space of differential operators.

Notice that ${\mathcal D}$ has a direct sum decomposition of the form
${\mathcal D}=\oplus_{n\geq 0} {\mathcal D}_n$
where ${\mathcal D}_n$ is the set of differential operators of order $n$ with no lower-order terms. Note that $\mathcal{D}_n$  is invariant under the map induced by $t_R$ and finitely generated over $R.$
So
$$E(k)=\bigoplus_{n\geq 0}E(k)_n$$ 
where $E(k)_n={\mathcal D}_m/\m {\mathcal D}_n$ and each $E(k)_n$ is  equivariant  and finitely generated $R-$module. This gives a decomposition as  follows:
$$H^i_{\m}(M)\simeq\bigoplus_{n\geq 0}\mathrm{Hom}_R(\mathrm{Ext}^{n-i}_R(M, R), E(k)_n)$$ 
But  each  factor  $\mathrm{Hom}_R(\mathrm{Ext}^{n-i}_R(M, R), E(k)_n)$ is finitely generated over $R$ so 
\ref{rema} applies to it.  Since there are countably many factors in the above decomposition, we see $\Spec(t^i)$ is countable. 
\hfill\endproof

\hfill

Now \ref{mainresult} is completely proven.

\subsection{Degeneration of the Dolbeault-Fr\"olicher spectral sequence
with coefficients in a local system}

The next result, interesting in itself, proves that  on compact LCK manifolds with proper potential, in the Dolbeault-Fr\"olicher spectral sequence with coefficients in a local system $L_\alpha$,
$$E_1^{p, q}:=H^q(M, \Omega_M^p\otimes L_\alpha)\Rightarrow H^{p+q}(M, L_\alpha(\C)),$$
all the terms vanish at $E_2$ level: $E_2^{p,q}=0$
(where $L_\al(\C)$ denotes the local system associated to $L_\al$). This parallels the degeneration of this spectral sequence at $E_1$ level for compact K\"ahler manifolds (where $L_\al$ is taken to be trivial). In particular, this gives a new proof to \ref{MN} and  produces new examples of compact complex manifolds that do not carry LCK metrics with potential.
One of the approaches to finding such manifolds is due
to S. Rollenske (\cite{Rol}), who showed that on a nilmanifold, the
Dolbeault-Fr\"olicher spectral sequence does not necessarily degenerate,
and gave examples when the $n$-th differential $d_n$ is
non-zero, for
arbitrarily high $n$.

\hfill

\proposition\label{_Poincare_weighted_Lemma_}
Let $M$ be a compact LCK manifold with proper potential,
 $\alpha\in \Hom(\Gamma, \R_+)$ a positive character,
and $L_\alpha$ the corresponding line bundle.
For any $p, q$, consider the map 
$$\6 _{p.q}:\; H^p(M, \Omega^q_M \otimes L_\alpha)\arrow H^p(M, \Omega^{q+1}_M \otimes L_\alpha).$$
Then $\ker \6 _{q, p+1}=\im\6 _{q,p}$, for all $p, q$.

\hfill

\noindent{\bf Proof:}
The monodromy map $\tilde t$ on $\tilde M$
is the exponential of a holomorphic vector field
$X$. This is proven in \cite[Theorem 2.3]{_OV:flow_} using the 
embedding of $M$ in a Hopf manifold $\C^N\setminus\{0\})/\langle A\rangle$ where $A$ is linear, with all eigenvalues smaller than 1. The holomorphic vector field is then $X=\log A$. In particular:
$$\tilde t^*(\eta)=\Lie_X\eta.$$
Let  now $[\eta]\in H^p(M, \Omega^{q+1}_M \otimes L_\alpha).$ A representative $\eta$ can be seen as a $(q+1, p)-$form on $\tilde M$ which is {$\overline{\6}$-closed}  and automorphic of weight $\alpha$.

Suppose   $\eta$  is also $\6$-closed. Then,  
since $\tilde t^*(\eta)=\alpha \cdot \eta$, we obtain 
$$\alpha \cdot \eta=\Lie_X(\eta)=di_X\eta+i_Xd\eta,$$
by Cartan's formula. 

But $\overline{\6}(\eta)=\6\eta=0$ by assumption,  thus $i_Xd\eta=0$, and we are left with:
$$\alpha \cdot \eta=\6(i_X\eta)+{\overline \6}(i_X\eta).$$
As $X$ is holomorphic, $i_X\eta$ is of type $(q,p)$, and hence ${\overline \6}(i_X\eta)$ is of type $(q, p+1).$ On the other hand both $\6 i_X(\eta)$  and  $\alpha \cdot \eta$ are of type $(q+1, p)$, implying  
${\overline \6}(i_X\eta)=0$ 
and 
$$\alpha \cdot \eta =\6 (i_X\eta).$$
This yields 
$\eta =\6 i_X\left(\frac{1}{\alpha}\eta\right)$, and hence  
$\eta\in \im(\6_{q,p}).$ \endproof

\section{Weighted Bott-Chern cohomology for LCK manifolds with potential}

We now generalize \cite[Theorem 4.7]{_OV:MN_}. We have:

\hfill

\proposition \label{_BC_via_Dolbeault_Proposition_}
Let $(M,I,g)$ be a compact LCK manifold. Then the following sequence is exact  for all
$\al\in \C$ but a discrete countable subset:
\begin{equation}\label{exgen}
 H^{q-1}_{\bar \6}(\Omega^p_M\otimes L_\al)\oplus \overline{H^{p-1}_{\bar \6}(\Omega^q_M\otimes L_\al)}\stackrel{\6_\theta+\bar\6_\theta}{\longrightarrow}H^{p,q}_{BC}(M,L_\al)\stackrel{\nu}{\longrightarrow}H^{p+q}(M,L_\al(\C))
\end{equation}
where $\nu$ is the tautological map, $\6_\theta=\6-\theta^{1,0}$ and $\bar\6_\theta=\bar\6-\theta^{0,1}$.

\hfill

\noindent{\bf Proof:} We prove that
$\im(\6_\theta+\bar\6_\theta)=\ker \nu$. Let $\eta$ be a
$(p,q)$-form with values in $L_\al$ whose class vanishes
in the cohomology of the local system $L_\al(\C)$. Then
$\eta=d_\theta\be$.
Suppose that $\beta$ has only two Hodge components,
$\beta=\be^{p,q-1}+\be^{p-1,q}$.
Then $\eta$ decomposes as
$\eta=\bar\6_\theta \be^{p,q-1}+\6_\theta\be^{p-1,q}$.
On the other hand, as $\eta$ is of bidegree $(p,q)$, we have
$\6_\theta \be^{p,q-1}=0$ and $\bar\6_\theta
\be^{p-1,q}=0$, and hence $\be^{p,q-1}$ and $\be^{p-1,q}$
produce the cohomology classes in
$[\be^{p-1,q}]\in H^{q-1}_{\bar \6}(\Omega^p_M\otimes L_\al)$ and
$[\be^{p,q-1}]\in \overline{H^{p-1}_{\bar    \6}(\Omega^q_M\otimes L_\al)}$.
Then
$[\eta]_{BC}=\6_\theta
[\be^{p-1,q}]+\bar\6_\theta[\be^{p,q-1}]$.

It remains to reduce \ref{_BC_via_Dolbeault_Proposition_}
to the case when $\beta$ has only two Hodge components. We
may already assume that $H^{p,q}(L_\alpha)=0$ for all $p
,q$ (\ref{mainresult}).
We use induction by the number of Hodge components.
Take the outermost Hodge component of
$\beta$, say, $\beta^{p-d-1, q+d}$, with $d> 0$. 
Then  $\bar\6_\theta(\beta^{p-d-1, q+d})=0$, hence, by vanishing
of the Dolbeault cohomology group $H^{p-d-1,q+d}(L_\alpha)$,
we have  $\beta^{p-d-1, q+d}=\bar\6_\theta(\gamma)$, where
$\gamma\in \Lambda^{p-d-1, q-1+d}(M, L_\alpha)$ is an
$L_\alpha$-valued $(p-d-1, q-1+d)$-form.
Now if we replace
$\beta$ by $\beta-d_\theta\gamma$, we obtain another form $\beta'$ such that
$\eta=d_\theta\beta'$, and $\beta'$ has a smaller number
of Hodge components. 
\endproof

\hfill 

As compact LCK manifolds with potential are topologically equivalent with Vaisman manifolds, \ref{pot} (ii), by \ref{MN} their cohomology of the local system $L_\al(\C)$ vanishes identically. Together with our main result (\ref{mainresult}), this proves the following generic vanishing of Bott-Chern cohomology (we keep the notations in Section 3):

\hfill

\corollary \label{_BC_vanishing_Corollary_}
Let $M$ be an LCK manifold with proper potential,
 $\alpha\in \C$ and $L_\al$ the flat line bundle
corresponding to $\al\cdot\theta$. Then
$H^{p,q}_{BC}(M,L_\al)=0$ for all
$\al\in \C$ but a discrete countable subset.

 \hfill

 \remark Note that $H^{p,q}_{BC}(M,L_\alpha)=0$
 implies the $d_{\al\theta}d_{\al\theta}^c$-lemma at the level $(p,q)$, and
hence our result 
says that, generically,  a compact LCK manifold with proper 
potential satisfies the $d_{\al\theta}d_{\al\theta}^c$-lemma for all $(p,q)$.

\hfill

\noindent{\bf Acknowledgement.} L.O. and V.V. are grateful to Higher School of Economics, Moscow, and M.V. is grateful to the University of Bucharest and ICUB for facilitating mutual visits during which parts of this research was carried on.

The authors thank L. Positselski for his help on Matlis duality, A. Otiman for pointing out an incomplete argument, and   the anonymous referee for a very careful reading of the manuscript and for her or his pertinent, useful remarks.

{\small

\noindent {\sc Liviu Ornea\\
University of Bucharest, Faculty of Mathematics, \\14
Academiei str., 70109 Bucharest, Romania}, and:\\
{\sc Institute of Mathematics "Simion Stoilow" of the Romanian
Academy,\\
21, Calea Grivitei Str.
010702-Bucharest, Romania\\
\tt lornea@fmi.unibuc.ro, \ \  Liviu.Ornea@imar.ro

\hfill

\noindent {\sc Misha Verbitsky\\
Laboratory of Algebraic Geometry, \\
Faculty of Mathematics, National Research University HSE,\\
7 Vavilova Str. Moscow, Russia}, also: \\
{\sc Universit\'e Libre de Bruxelles, D\'epartement de Math\'ematique\\
Campus de la Plaine, C.P. 218/01, Boulevard du Triomphe\\
B-1050 Brussels, Belgium\\
\tt verbit@verbit.ru }

\hfill

\noindent {\sc Victor Vuletescu\\
University of Bucharest, Faculty of Mathematics, \\14
Academiei str., 70109 Bucharest, Romania.}\\
\tt vuli@fmi.unibuc.ro
}}


\end{document}